\documentstyle[12pt,subeqn]{article}
\newfont{\twelvemsb}{msbm10 scaled\magstep1}
\newfont{\eightmsb}{msbm8}
\newfam\msbfam
\textfont\msbfam=\twelvemsb
\scriptfont\msbfam=\eightmsb
\catcode`\@=11
\def\Bbb{\ifmmode\let\next\Bbb@\else
  \def\next{\errmessage{Use \string\Bbb\space only in math mode}}\fi\next}
\def\Bbb@#1{{\fam\msbfam{{#1}}}}

\newcommand{\sect}[1]{\setcounter{equation}{0}\section{#1}}

\newcommand{\be}{\begin{equation}}
\newcommand{\ee}{\end{equation}}
\newcommand{\bea}{\begin{eqnarray}}
\newcommand{\eea}{\end{eqnarray}}
\newcommand{\beo}{\begin{eqnarray*}}
\newcommand{\eeo}{\end{eqnarray*}}

\newcommand{\cA}{{\cal A}}
\newcommand{\cM}{{\cal M}}
\newcommand{\cR}{{\cal R}}
\newcommand{\cY}{{\cal Y}}
\newcommand{\cZ}{{\cal Z}}
\newcommand{\CC}{{\Bbb C}}
\newcommand{\II}{{\Bbb I}}
\newcommand{\ZZ}{{\Bbb Z}}
\newcommand{\eps}{{\varepsilon}}
\newcommand{\elp}{{\cA_{q,p}(\widehat{sl}(2)_c)}}
\newcommand{\snh}{{\mbox{snh}}}
\newcommand{\tr}{{\mbox{Tr}}}

\newtheorem{defin}{Definition}
\newtheorem{prop}{Proposition}
\newtheorem{thm}{Theorem}

\begin{document}
\newpage
\pagestyle{empty}
\setcounter{page}{0}
\vfill
\begin{center}

{\Large {\bf Poisson structures on the center 

 \vspace{5mm}

of the elliptic algebra $\elp$}}

\vspace{7mm}

{\large J. Avan}

\vspace{4mm}

{\em LPTHE, CNRS-URA 280, Universit\'es Paris VI/VII, France}

\vspace{7mm}

{\large L. Frappat \footnote{On leave of absence from Laboratoire de 
Physique Th\'eorique  ENSLAPP.}}

\vspace{4mm}

{\em Centre de Recherches Math\'ematiques, Universit\'e de Montr\'eal, Canada}

\vspace{7mm}

{\large M. Rossi, P. Sorba}

\vspace{4mm}

{\em Laboratoire de Physique Th\'eorique ENSLAPP \footnote{URA 1436 du
CNRS, associ\'ee \`a l'\'Ecole Normale Sup\'erieure de Lyon et \`a
l'Universit\'e de Savoie.}\\
Annecy-le-Vieux C\'edex 
and ENS Lyon, Lyon C\'edex 07, France}

\end{center}

\vfill

\begin{abstract}
It is shown that the elliptic algebra $\elp$ has a non trivial center
at the critical level $c=-2$, generalizing the result of Reshetikhin
and Semenov-Tian-Shansky for trigonometric algebras. A family of Poisson
structures indexed by a non-negative integer $k$ is constructed on this center.
\end{abstract}

\def\abstractname{R\'esum\'e}
\begin{abstract}
On montre que l'alg\`ebre elliptique $\elp$ a un centre non trivial au
niveau critique $c=-2$, g\'en\'eralisant le r\'esultat de Reshetikhin
et Semenov-Tian-Shansky pour les alg\`ebres trigonom\'etriques. On
construit sur ce centre une famille de structures de Poisson indic\'ees 
par un entier positif $k$.
\end{abstract}

\vfill
\vfill

\rightline{ENSLAPP-AL-644/97}
\rightline{CRM-2474}
\rightline{PAR-LPTHE 97-17}
\rightline{q-alg/9705012}
\rightline{May 1997}

\newpage
\pagestyle{plain}

\sect{Introduction}

The concept of deformed Virasoro and $W_N$ algebras has recently arisen in
connection with several aspects of integrable systems, both in quantum and
statistical mechanics. On the one hand, $q$-deformed Virasoro and $W_N$
algebras of operators were introduced \cite{SKAO,AKOS} as an extension of
the Virasoro and $W_N$ algebras identified in the quantum Calogero-Moser
model \cite{AMOS}. In the same way as Jack polynomials (eigenfunctions of
the quantum Calogero-Moser model) arise as singular vectors of $W_N$ 
algebras \cite{AMOS,MaChe}, Mac-Donald polynomials (eigenfunctions of
the relativistic Ruijsenaars-Schneider model) are singular vectors of the
$q$-deformed $W_N$ algebras. The parameter $q$ is connected to the 
supplementary parameter (speed of light) in Ruijsenaars-Schneider models 
accounting for relativistic invariance.

On the other hand, it was shown that $q$-deformed Virasoro algebras provided
a dynamical symmetry for the Andrew-Baxter-Forrester restriction of RSOS
models in statistical mechanics \cite{LP}. This model is characterized by a 
matrix of Boltzmann weights which is an elliptic solution of the Yang-Baxter
equation \cite{ABF} (see also \cite{FJMOP} for further developments).

These algebras were shown to arise in fact from a systematic procedure of
construction mimicking the already known scheme \cite{FF1} for Virasoro
and $W_N$ algebras: undeformed algebras can be obtained by quantization of
a Poisson structure on the center of the enveloping algebra 
$U(\widehat{sl}(N)_c)$ of an affine algebra $\widehat{sl}(N)_c$ at the
critical level $c=-N$ \cite{FF1}. This construction is easier to achieve
using the bosonized representation \cite{Wak} of $U(\widehat{sl}(N)_c)$,
leading to the well-known Miura transformation formula for Virasoro and
$W_N$ generators.

The procedure for $q$-deformed Virasoro and $W_N$ algebras is similar. One
starts from the quantum affine Lie algebra $U_q(\widehat{sl}(N)_c)$ 
\cite{Dr,Jim}. At the critical value $c=-N$, the algebra has a 
multidimensional center \cite{RSTS} on which a Poisson structure can be
defined as limit of the commutator structure \cite{FR}. This Poisson
structure is the semi-classical limit of the $q$-deformed Virasoro (for
$N=2$) or $W_N$ algebras. Its quantization \cite{FF2} gives rise to these
$q$-deformed algebras. In fact, the construction of classical and quantized
algebras was not achieved in \cite{FF2} by direct computation but using the
$q$-deformed bozonization \cite{AOS} of $U_q(\widehat{sl}(N)_c)$ and the
corresponding $q$-Miura transformation. Interestingly the $q$-deformed $W_N$ 
algebras are characterized by elliptic structure coefficients, a fact on 
which we shall comment at the end: for instance, the $q$-deformed Virasoro
algebra is defined by the generating operator $T(z)$ such that
\be
f_{1,2}(w/z) \, T(z) \, T(w) - f_{1,2}(z/w) \, T(w) \, T(z) 
= \frac{(1-q)(1-p/q)}{1-p} \left( \delta\Big(\frac{w}{zp}\Big) 
- \delta\Big(\frac{wp}{z}\Big) \right)\, ,
\ee
where
\be
f_{1,2}(x) = \frac{1}{1-x} ~ \frac{(x|q,pq^{-1};p^2)_\infty}
{(x|pq,p^2q^{-1};p^2)_\infty}
\ee
and as usual
\be
(x|a_1 \dots a_k;t)_\infty \equiv \prod_{i=1}^k \prod_{n=0}^\infty (1-a_ixt^n)
\, .
\ee

The parameters $p$ and $q$ are rewritten as $q=e^h$, $p=e^{h(1-\beta)}$; 
hence $h$ is the deformation parameter and $\beta$ is the ``quantization''
parameter, the semi-calssical limit $\beta \rightarrow 0$ giving back the 
$q$-deformed Poisson bracket and the limit $h \rightarrow 0$ giving back the 
linear Virasoro algebra.

It is a natural question to extend these constructions to other deformations
of affine Lie algebras, in particular to the so-called elliptic quantum
algebras introduced by \cite{FIJKMY} and further studied in \cite{KLP,JMK}.
Defined in a way similar to \cite{RSTS}, they use a structure $R$-matrix with
elliptic dependence associated to the eight-vertex model \cite{Ba}.

As a first step in this extension, our purpose here is to examine the classical
Poisson structure on the center of the elliptic quantum algebra $\elp$. 
In a first step we establish the existence of a non-trivial center for this 
algebra when $c=-2$ and derive the relevant exchange algebra relations. In a 
second step we examine the limit $c \rightarrow -2$ and derive the Poisson 
structure directly as a limit of this exchange algebra. The direct method is 
here necessary since no bosonized version \`a la Wakimoto is get available for 
this algebra.

\sect{The elliptic quantum algebra $\elp$}

{\bf \thesection.1.} 
Consider the $R$-matrix of the eight vertex model found by Baxter
\cite{Ba}:
\be
R_{12}(x) = \frac{1}{\mu(x)} \left(\begin{array}{cccc}
a(u) & 0 & 0 & d(u) \cr 0 & b(u) & c(u) & 0 \cr 
0 & c(u) & b(u) & 0 \cr d(u) & 0 & 0 & a(u) \cr 
\end{array}\right) \label{eq21}
\ee
where the functions $a(u), b(u), c(u), d(u)$ are given by 
\be
a(u) = \frac{\snh(\lambda-u)}{\snh(\lambda)} \,, \quad
b(u) = \frac{\snh(u)}{\snh(\lambda)} \,, \quad
c(u) = 1 \,, \quad
d(u) = k\,\snh(\lambda-u)\snh(u) \,.
\ee
The function $\snh(u)$ is defined by $\snh(u) = -i\mbox{sn}(iu)$ where 
$\mbox{sn}(u)$ is Jacobi's elliptic function with modulus $k$.
If the elliptic integrals are denoted by $K,K'$ (let ${k'}^2=1-k^2$),
\be
K = \int_0^1 \frac{dx}{\sqrt{(1-x^2)(1-k^2x^2)}} \qquad\mbox{and}\qquad
K' = \int_0^1 \frac{dx}{\sqrt{(1-x^2)(1-{k'}^2x^2)}} \,,
\ee
the functions $a(u), b(u), c(u), d(u)$ can be seen as functions of the
variables
\be
p = \exp\Big(-\frac{\pi K'}{K}\Big) \,, \qquad
q = - \exp\Big(-\frac{\pi\lambda}{2K}\Big) \,, \qquad
x = \exp\Big(\frac{\pi u}{2K}\Big) \,.
\ee
The normalization factor $\mu(x)$ is chosen as follows \cite{JMK}:
\bea
&& \frac{1}{\mu(x)} = \frac{1}{\kappa(x^2)} \frac{(p^2;p^2)_\infty}
{(p;p)_\infty^2} \frac{\Theta_{p^2}(px^2)\Theta_{p^2}(q^2)}
{\Theta_{p^2}(q^2x^2)} \,, \\
&& \frac{1}{\kappa(x^2)} = \frac{(q^4x^{-2};p,q^4)_\infty
(q^2x^2;p,q^4)_\infty (px^{-2};p,q^4)_\infty (pq^2x^2;p,q^4)_\infty}
{(q^4x^2;p,q^4)_\infty (q^2x^{-2};p,q^4)_\infty (px^2;p,q^4)_\infty 
(pq^2x^{-2};p,q^4)_\infty} \,,
\eea
where one defines the infinite multiple products as usual by
\be
(x;p_1,\dots,p_m)_\infty = \prod_{n_i \ge 0} (1-xp_1^{n_1} \dots p_m^{n_m})
\ee
and $\Theta$ is the Jacobi $\Theta$ function:
\be
\Theta_{p^2}(x) = (x;p^2)_\infty (p^2x^{-1};p^2)_\infty (p^2;p^2)_\infty \,.
\ee
\begin{prop}\label{propone}
The matrix $R_{12}$, element of $End(\CC^2) \otimes End(\CC^2)$, has the 
following properties:
\\
\indent -- unitarity: $R_{21}(x^{-1}) R_{12}(x) = 1$,
\\
\indent -- crossing symmetry: $R_{21}(x^{-1})^{t_1} 
= (\sigma^1 \otimes \II) R_{12}(-q^{-1}x) (\sigma^1 \otimes \II)$,
\\
\indent -- antisymmetry: $R_{12}(-x) = - (\sigma^3 \otimes \II) 
R_{12}(x) (\sigma^3 \otimes \II)$,
\\
where $\sigma^1,\sigma^2,\sigma^3$ are the $2 \times 2$ Pauli matrices and
$t_i$ denotes the transposition in the space $i$.
\end{prop}
The proof is straightforward by direct calculation.

\medskip

For the definition of the elliptic quantum algebra $\elp$, we need to
use a modified $R$-matrix $R_{12}^+(x)$ defined by 
\be
R_{12}^+(x) = \tau(q^{1/2}x^{-1}) R_{12}(x) \, ,
\ee
where the prefactor $\tau$ is given by
\be
\tau(x) = x^{-1} \frac{(qx^2;q^4)_\infty (q^3x^{-2};q^4)_\infty}
{(qx^{-2};q^4)_\infty (q^3x^2;q^4)_\infty} \,.
\ee
Let us stress that the function $\tau$ is periodic: $\tau(x) = \tau(xq^2)$ and
$R_{12}^+(x)$ obeys a quasi-periodicity property:
\be
R_{12}^+(-p^{\frac {1}{2}}x)= (\sigma^1 \otimes \II) 
\left (R_{21}^+(x^{-1})\right )^{-1} (\sigma^1 \otimes \II) \, .
\ee
The crossing symmetry and the unitarity property of $R_{12}$ (see Proposition
\ref{propone}) then allow exchange of inversion and transposition for the
matrix $R_{12}^+$ as:
\be
\Big(R_{12}^+(x)^{t_2}\Big)^{-1} = \Big(R_{12}^+(q^2x)^{-1}\Big)^{t_2}\, .
\label{eq211}
\ee

\medskip

{\bf \thesection.2.} 
The elliptic quantum algebra $\elp$ has been introduced by \cite{FIJKMY,JMK}.
The elliptic quantum algebra $\elp$ is an algebra of operators 
$L_{\eps\eps',n}$ such that ($\eps,\eps'=+$ or $-$)
\subequations
\bea
&& \mbox{i)} ~ L_{\eps\eps',n} = 0 \mbox{ if } \eps\eps' \ne (-1)^n \,, 
\label{eq212a} \\
&& \mbox{ii) one defines } L_{\eps\eps'}(z) = \sum_{n\in\ZZ} 
L_{\eps\eps',n} z^n \,, \label{eq212b}
\eea
\endsubequations
and encapsulate them into a $2 \times 2$ matrix
\be
L = \left(\begin{array}{cc} L_{++} & L_{+-} \cr L_{-+} & L_{--} \cr
\end{array}\right) \,.
\ee
Note that by virtue of (\ref{eq212a}-b), $L_{++}(z)$ and $L_{--}(z)$ are even
while $L_{+-}(z)$ and $L_{-+}(z)$ are odd functions of $z$.

\medskip

One then defines $\cA_{q,p}(\widehat{gl}(2)_c)$ by imposing the following 
constraints on $L_{\eps\eps'}(z)$:
\be
R_{12}^+(z/w) \, L_1(z) \, L_2(w) = L_2(w) \, L_1(z) \, R_{12}^{+*}(z/w) \,,
\label{eq214}
\ee
where $L_1(z) \equiv L(z) \otimes \II$, $L_2(z) \equiv \II \otimes L(z)$ and
$R^{+*}_{12}$ is defined by $R^{+*}_{12}(x,q,p) \equiv 
R^+_{12}(x,q,pq^{-2c})$. This shift plays an essential role in establishing 
the existence of a center at $c=-2$.
\\
Since the $q$-determinant $q$-$\det L(z) \equiv L_{++}(q^{-1}z) L_{--}(z) - 
L_{-+}(q^{-1}z) L_{+-}(z)$ is in the center of $\cA_{q,p}(\widehat{gl}(2)_c)$,
it can be ``factored out'', being set to the value $q^{c/2}$ so as to get
$\elp = \cA_{q,p}(\widehat{gl}(2)_c)/
\langle q\mbox{-}\det L - q^{c/2} \rangle$.

It is useful, both for our next computations and the establishing of the 
$p \rightarrow 0$ trigonometric limit, to introduce the following two 
matrices:
\be
L^+(z) \equiv L(q^{c/2}z) \,, \qquad
L^-(z) \equiv \sigma^1 L(-p^{1/2}z) \sigma^1 \,.
\ee
They obey coupled constraint relations following from (\ref{eq214}) and 
periodicity/unitarity properties of the matrix $R^+_{12}$:
\be
\begin{array}{l}
R^+_{12}(z/w) ~ L^\pm_1(z) ~ L^\pm_2(w) = 
L^\pm_2(w) ~ L^\pm_1(z) ~ R^{+*}_{12}(z/w) \,, \\
\bigg. R^+_{12}(q^{c/2}z/w) ~ L^+_1(z) ~ L^-_2(w) = 
L^-_2(w) ~ L^+_1(z) ~ R^{+*}_{12}(q^{-c/2}z/w) \,.
\end{array} \label{eq216}
\ee
Note that a further renormalization of the modes 
$\bar L_n=\left (-p^{-\frac {1}{2}} \right ) ^{\max (-n,0)}L_n$ is required 
in order to get the (degenerate) trigonometric algebra 
$U_q(\widehat{sl}(2)_c)$ from $\elp$ when $p\rightarrow 0$. It ensures the 
vanishing of half of the degrees of freedom in $L^\pm(z)$, but decouples 
completely $L^+(z)$ from $L^-(z)$, thereby keeping the same overall 
dimensionality \cite{JMK}. It is therefore a non-trivial, ``discontinuous'' 
procedure.

We now state the first important result of our study.

\sect{The center of the elliptic quantum algebra $\elp$}

\begin{thm}\label{thmone}
At $c=-2$, the operators generated by 
\be
t(z) = {\mbox Tr}(L(z)) = {\mbox Tr}\Big(L^+(q^{c/2}z) L^-(z)^{-1}\Big) \label{eq31}
\ee
commute with the algebra $\elp$.
\end{thm}
The formula for $t(z)$ in the elliptic case is exactly identical to the one in 
the trigonometric case \cite{RSTS}. The proof follows on similar lines and is 
interesting to work out in more detail since it involves all properties of the 
$R$-matrix and uses tricks to be similarly used later.

\medskip

\noindent
{\bf Proof:}

$\bullet$ step 1: One derives further useful exchange relations from 
eq. (\ref{eq216}). 
\\
Defining $\widetilde{L}^\pm(z) \equiv (L^\pm(z)^{-1})^t$
(this definition is unambiguous at least for such $2 \times 2$ 
matrices of operators, as can be seen using Borel decomposition), one has:
\be
(R_{12}^+(z/w)^{t_2})^{-1} \, L_1^\pm(z) \, \widetilde{L}_2^\pm(w) = 
\widetilde{L}_2^\pm(w) \, L_1^\pm(z) \, (R_{12}^{+*}(z/w)^{t_2})^{-1}\, .
\label{eq32}
\ee
This follows from starting with $R_{12}^+  L_1^\pm  L_2^\pm = 
L_2^\pm L_1^\pm  R_{12}^{+*}$, rewriting it as
$L_1^\pm  R_{12}^{+*}  (L_2^\pm)^{-1} = (L_2^\pm)^{-1} 
R_{12}^{+*}  L_1^\pm$, transposing with respect to indices 2 (which is 
allowed since $R_{12}$ is a $c$-number matrix) to obtain
$L_1^\pm \, ((L_2^\pm)^{-1})^{t_2} \, (R_{12}^{+*})^{t_2} = 
(R_{12}^{+*})^{t_2} \, ((L_2^\pm)^{-1})^{t_2} \, L_1^\pm$
and finally getting (\ref{eq32}) by multiplying l.h.s. and r.h.s. by suitable 
inverses of $(R_{12}^+)^{t_2}$.
\\
Similarly, one gets
\be
(R_{12}^+(q^{c/2}z/w)^{t_2})^{-1} \, L_1^+(z) \, \widetilde{L}_2^-(w) = 
\widetilde{L}_2^-(w) \, L_1^+(z) \, (R_{12}^{+*}(q^{-c/2}z/w)^{t_2})^{-1}\, ,
\label{eq33}
\ee
by applying a similar scheme of rewritings to eq. (\ref{eq216}).

\medskip

$\bullet$ step 2: To compute $\Big[ t(z),L^+(w) \Big]$, one rewrites at
$c=-2$
\bea
t(z) \, L_2^+(w) &=& \tr_1\Big( L_1^+(zq^{-1}) \widetilde{L}_1^-(z)^{t_1} \Big) 
L_2^+(w) \nonumber \\
&=& \tr_1 \Big( L_1^+(zq^{-1})^{t_1} \widetilde{L}_1^-(z) L_2^+(w) \Big)\, ,
\eea
since one is allowed to exchange transposition under a trace procedure.
One then commutes $L_2^+(w)$ through $\widetilde{L}_1^-(z)$ using eq.
(\ref{eq33}) as:
\be
t(z) \, L_2^+(w) = \tr_1\Big( L_1^+(zq^{-1})^{t_1} (R_{21}^+(w/qz)^{t_1})^{-1} 
L_2^+(w) \widetilde{L}_1^-(z) R_{21}^{+*}(qw/z)^{t_1} \Big)\, .
\label{eq35}
\ee

\medskip

$\bullet$ step 3: One uses the exchange algebra (\ref{eq216}) to get a suitable
expression for the first three terms in (\ref{eq35}):
\\
- transposition with respect to indices 1 gives:
$L_1^+(z)^{t_1}  R_{12}^+(z/w)^{t_1}  L_2^+(w) = L_2^+(w) 
 {R_{12}^{+*}(z/w)}^{t_1}  L_1^+(z)^{t_1}$;
\\
- unitarity property and redefinition of $z$ as $zq^{-1}$ give:
\be
L_1^+(zq^{-1})^{t_1} \, (R_{21}^+(qw/z)^{-1})^{t_1} \, L_2^+(w) = L_2^+(w) 
\, (R_{21}^{+*}(qw/z)^{-1})^{t_1} \, L_1^+(zq^{-1})^{t_1} \,; \label{eq36}
\ee
Eq. (\ref{eq211}) inserted in (\ref{eq36}) gives (after due exchange of 
spaces 1 and 2):
\be
L_1^+(zq^{-1})^{t_1} \, (R_{21}^+(w/qz)^{t_1})^{-1} \, L_2^+(w) = L_2^+(w) 
\, ({R_{21}^{+*}(w/qz)}^{t_1})^{-1} \, L_1^+(zq^{-1})^{t_1} \,. \label{eq37}
\ee

\medskip

$\bullet$ step 4: One now inserts (\ref{eq37}) into (\ref{eq35}) to get
\bea
t(z) \, L_2^+(w) &=& L_2^+(w) \, \tr_1\Big( (R_{21}^{+*}(w/qz)^{t_1})^{-1} 
L_1^+(zq^{-1})^{t_1} \widetilde{L}_1^-(z) R_{21}^{+*}(qw/z)^{t_1} \Big)
\nonumber \\
&=& L_2^+(w) \, \tr_1\Big( (R_{21}^{+*}(qw/z)^{-1})^{t_1} 
L_1^+(zq^{-1})^{t_1} \widetilde{L}_1^-(z) R_{21}^{+*}(qw/z)^{t_1} \Big)\, .
\eea
One now needs to use the fact that under a trace over the space 1 one has
$\tr_1 \Big( R_{21} Q_1 {R'}_{21} \Big) = \tr_1 \Big( Q_1 {R'_{21}}^{t_2} 
{R_{21}}^{t_2} \Big)^{t_2}$:
\be
t(z) \, L_2^+(w) = L_2^+(w) \, \tr_1\Big( L_1^+(zq^{-1})^{t_1} 
\widetilde{L}_1^-(z) R_{21}^{+*}(wq/z)^{t_1t_2} 
(R_{21}^{+*}(wq/z)^{-1})^{t_1t_2} \Big)^{t_2}\, .
\ee
The last two terms in the right hand side cancel each other, leaving a trivial
dependence in space 2 and 
$\tr_1\Big( {L_1^+(zq^{-1})}^{t_1} \widetilde{L}_1^-(z)\Big) \equiv t(z)$ 
in space 1. This shows the commutation of $t(z)$
with $L^+(w)$ and therefore with $L^-(w)=\sigma ^1 L^+(-p^{\frac {1}{2}}q^{-\frac {c}{2}}w) \sigma ^1$, hence with the full algebra $\elp$ at $c=-2$.
\hfill \rule{5pt}{5pt}

\medskip

It is very important to notice here that owing to the parity properties of 
the modes of $L(z)$ defined above and the identical properties of the modes of
$L(z)^{-1}$ or $\widetilde{L}(z)$, the operator $t(z)$ is even in $z$ and has
therefore only {\em even} modes.

\medskip

We now study the specific behaviour of the exchange algebra of $t(z)$ with
$t(w)$ in the neighborhood of $c=-2$. Notice that we have not proved that 
$t(z)$ exhausted the center of $\elp$ at $c=-2$. Hence the exchange algebra
may not close on $t(z),t(w)$ in the neighborhood of $c=-2$. This will however 
turn out to be true.

\sect{Poisson algebra of $t(z)$}

By virtue of Theorem \ref{thmone} and by the definition (\ref{eq31}) of 
$t(z)$ itself, the elements $t(z)$ and $t(w)$ are mutually commuting at the
critical level $c=-2$. This implies a natural Poisson structure on the
center $\cZ$: if $\Big[ t(z),t(w) \Big] = (c+2) \ell(z,w) + o(c+2)$, then
a Poisson bracket on $\cZ$ is yielded by 
$\Big\{ t(z)_{cr},t(w)_{cr} \Big\} = \ell(z,w)_{cr}$
(the subscript ``cr'' meaning that all expressions are taken at $c=-2$).
We now state the main result of the paper:
\begin{thm}\label{thmtwo}
The elements $t(z)$ form a closed algebra under the natural Poisson 
bracket on $\cZ$. More precisely, we have (we suppress the subscript ``$cr$'' 
for simplicity)
\be
\Big\{ t(z),t(w) \Big\} = -(\ln q) \left((w/z)\frac{d}{d(w/z)} 
\ln\tau(q^{1/2}w/z) - (z/w)\frac{d}{d(z/w)}\ln\tau(q^{1/2}z/w)\right) 
~ t(z)t(w) \,. \label{eq41}
\ee
\end{thm}
{\bf Proof:}

$\bullet$ step 1: one computes the exchange algebra between the operators
$t(z)$ and $t(w)$. From the definition of the element $t(z)$, one has
\be
t(z)t(w) = L(z)^{i_1}_{i_1} \, L(w)^{i_2}_{i_2} = 
L^+(z)^{j_1}_{i_1} \tilde L^-(z)^{j_1}_{i_1} \, 
L^+(w)^{j_2}_{i_2} \tilde L^-(w)^{j_2}_{i_2} \,.
\ee
The same kind of rewritings of eq. (\ref{eq216}) done in the first
step of the proof of Theorem \ref{thmone} leads to the following exchange
relation between the operators $\widetilde{L}^-$:
\be
R_{12}^+(z/w)^{t_1t_2} \, \widetilde{L}_1^-(z) \, \widetilde{L}_2^-(w) = 
\widetilde{L}_2^-(w) \, \widetilde{L}_1^-(z) \, R_{12}^{+*}(z/w)^{t_1t_2}
\label{eq43}\, .
\ee
The exchange relations (\ref{eq216}), (\ref{eq33}) and (\ref{eq43})
and the properties of the matrix $R_{12}$ given in proposition
\ref{propone} then allow us to move the matrices $L^+(w)$, $\tilde
L^-(w)$ to the left of the matrices $L^+(z)$, $\tilde L^-(z)$. One obtains
\be
t(z)t(w) = \cY(z/w)^{i_1i_2}_{j_1j_2} ~ L(w)^{j_2}_{i_2} ~ L(z)^{j_1}_{i_1} \,,
\label{eq44}
\ee
where the matrix $\cY(z/w)$ is given in terms of the matrix $R^+_{12}$ by
\be
\cY(z/w) = \bigg(\Big(R^+_{12}(w/z) \, {R^+_{12}(q^{c+2}w/z)}^{-1} \, 
{R^+_{12}(z/w)}^{-1}\Big)^{t_2} \, {R^+_{12}(q^cz/w)}^{t_2}\bigg)^{t_2} \,.
\label{eq45}
\ee
Making explicit the dependence of the matrix $\cY(z/w)$ on the function 
$\tau$, one gets
\be
\cY(z/w) = T(z/w) \cR(z/w) \,, \label{eq46}
\ee
where the matrix factor $\cR(z/w)$ depends only on the matrix (\ref{eq21}):
\be
\cR(z/w) = \bigg(\Big(R_{12}(w/z)\, R_{12}(z/q^{-c-2}w)
\, R_{12}(w/z)\Big)^{t_2} \, {R_{12}(z/q^cw)}^{t_2}\bigg)^{t_2}
\ee
and the numerical prefactor $T(z/w)$ contains all the $\tau$ dependence:
\be
T(z/w) = \frac{\tau(z/q^{1/2}w)\tau(w/q^{-c+1/2}z)}
{\tau(z/q^{-c-3/2}w)\tau(w/q^{1/2}z)} \,.
\ee
One easily checks the nice behaviour of $T(z/w)$ and $\cR(z/w)$ 
at $c=-2$:
\be
\begin{array}{l}
T(z/w)_{cr} = 1 \cr \bigg. \cR(z/w)_{cr} = \II_2 \otimes \II_2 \cr
\end{array} \qquad \Longrightarrow \qquad 
\cY(z/w)_{cr} = \II_2 \otimes \II_2 \,, \label{eq49}
\ee
making obvious the interest of the factorization (\ref{eq46}), so that one 
recovers $t(z)t(w) = t(w)t(z)$ at the critical level.

\medskip

$\bullet$ step 2: One computes the Poisson structure from the exchange
algebra (\ref{eq44}) in the neighborhood of $c=-2$.
As stressed at the beginning of the section, there is a natural Poisson
bracket on $\cZ$. From eqs. (\ref{eq44}) and (\ref{eq49}), one writes
\be
t(z)t(w) = t(w)t(z) + (c+2)\left(\frac{d\cY}{dc}(z/w)\right)^{i_1i_2}_{j_1j_2}
~ L(w)^{j_2}_{i_2} ~ L(z)^{j_1}_{i_1}\bigg\vert_{cr} + o(c+2)
\ee
and therefore
\be
\Big\{ t(z),t(w) \Big\} = \left(\frac{d\cY}{dc}(z/w)\right)^{i_1i_2}_{j_1j_2}
~ L(w)^{j_2}_{i_2} ~ L(z)^{j_1}_{i_1}\bigg\vert_{cr} \, . \label{eq411}
\ee
The main difficulty is of course the calculation of the derivative in 
(\ref{eq411}).
The equations (\ref{eq46}) and (\ref{eq49}) imply
\be
\frac{d\cY}{dc}(x)\bigg\vert_{cr} = \frac{dT}{dc}(x)\bigg\vert_{cr} 
\II_2 \otimes \II_2 ~+~ \frac{d\cR}{dc}(x)\bigg\vert_{cr} \,. \label{eq412}
\ee
Let us write $\cR(x)$ as 
\be
\cR(x) = \frac{1}{\mu(xq^c)\mu(xq^{-c-2})\mu(x^{-1})^2} ~ \cM(x) \,,
\ee
so that the dependence on the function $\mu$ is now explicit (note that 
$\cM(x)_{cr} = \mu(xq^{-2})\mu(x)\mu(x^{-1})^2$ by virtue of (\ref{eq49})). 
The matrix $\cM(x)$ is given by
\be
\cM(x) = \left(\begin{array}{cccc}
m_{11} & 0 & 0 & m_{22} \cr 0 & m_{12} & m_{21} & 0 \cr
0 & m_{21} & m_{12} & 0 \cr m_{22} & 0 & 0 & m_{11} \cr
\end{array}\right) \,,
\ee
the entries of which depending only on the elliptic functions 
$a(u), b(u), c(u), d(u)$ by
\subequations
\bea
m_{11} &=& a(xq^c)a(xq^{-c-2})a(x^{-1})^2 
+ 2a(xq^c)a(x^{-1})d(xq^{-c-2})d(x^{-1}) + a(xq^c)a(xq^{-c-2})d(x^{-1})^2 
\nonumber \\
&& - 2b(xq^{-c-2})b(x) + b(x)^2 + 1 \\
m_{22} &=& a(xq^{-c-2})a(x^{-1})^2 + 2a(x^{-1})d(xq^{-c-2})d(x^{-1})
+ a(xq^{-c-2})d(x^{-1})^2 \nonumber \\
&&- 2a(xq^c)b(xq^{-c-2})b(x) + a(xq^c)b(x)^2 + a(xq^c) \\
m_{12} &=& a(x^{-1})^2d(xq^c)d(xq^{-c-2}) 
+ 2a(xq^{-c-2})a(x^{-1})d(xq^c)d(x^{-1}) + b(xq^c)b(xq^{-c-2})b(x)^2 
\nonumber \\
&& - 2b(xq^c)b(x) + b(xq^c)b(xq^{-c-2}) + d(xq^c)d(xq^{-c-2})d(x^{-1})^2 \\
m_{21} &=& a(x^{-1})^2b(xq^c)d(xq^{-c-2}) 
+ 2a(xq^{-c-2})a(x^{-1})b(xq^c)d(x^{-1}) + b(xq^{-c-2})b(x)^2d(xq^c) 
\nonumber \\
&& - 2b(x)d(xq^c) + b(xq^{-c-2})d(xq^c) + b(xq^c)d(xq^{-c-2})d(x^{-1})^2 
\eea
\endsubequations
Using the formula
\be
\frac{d}{dc} f(xq^{2\alpha\pm(c+2)})\bigg\vert_{cr}  = \pm(\ln q) x\frac{d}{dx} f(xq^{2\alpha})
= \pm(2iK-\lambda) \frac{d}{du} f(u-2\alpha\lambda) \,, \label{eq415}
\ee
one shows after a long calculation, using various tricks in elliptic 
functions theory, that at $c=-2$ the derivatives of the diagonal entries of 
$\cM$ are equal while the derivatives of the off-diagonal entries 
identically vanish, namely:
\bea
&& \hspace{-5mm} \frac{d}{dc}m_{11}\bigg\vert_{cr} 
= \frac{d}{dc}m_{12}\bigg\vert_{cr} = 
(2iK-\lambda) \bigg((1-b(u+\lambda)^2)\frac{d}{du}b(u)^2 
- (1-b(u)^2)\frac{d}{du}b(u+\lambda)^2\bigg) \,, \\
&& \hspace{-5mm} \frac{d}{dc}m_{21}\bigg\vert_{cr} 
= \frac{d}{dc}m_{22}\bigg\vert_{cr} = 0 \,.
\eea
In particular the following relations are essential in this derivation:
\beo
&& \frac{\snh(a-u)\snh(a-v)-\snh(u)\snh(v)}{1-k^2\snh(u)\snh(v)\snh(a-u)
\snh(a-v)} = \snh(a-u-v)\snh(a) \,, \\
&& \frac{\snh(u)\snh(a-u)-\snh(v)\snh(a-v)}{\snh(u)\snh(a-v)-\snh(v)\snh(a-u)}
 = \frac{\snh(a-u-v)}{\snh(a)} \,, \\
&& \frac{\snh(u)\snh'(v)-\snh(v)\snh'(u)}{\snh^2(u)-\snh^2(v)} 
= \frac{1}{\snh(u+v)} \,.
\eeo
Moreover, one has $\mu(x)\mu(x^{-1}) = 1-b(u)^2$ and $\displaystyle{
\frac{\mu(x)}{\mu(xq^{-2})} = \frac{1-b(u)^2}{1-b(u+\lambda)^2}}$, which
imply
\bea
\frac{d}{dc} \left(\frac{1}{\mu(xq^c)\mu(xq^{-c-2})\mu(x^{-1})^2}\right) 
\bigg\vert_{cr} &=& (2iK-\lambda) \frac{1}{\mu(x)^2\mu(x^{-1})^2} 
\frac{d}{du} \left(\frac{\mu(x)}{\mu(xq^{-2})}\right) 
\nonumber \\
&=& -(2iK-\lambda) \frac{(1-b(u+\lambda)^2)\frac{d}{du}b(u)^2 
- (1-b(u)^2)\frac{d}{du}b(u+\lambda)^2}{(1-b(u)^2)^2 ~ (1-b(u+\lambda)^2)^2} 
\,. \nonumber \\
\eea
Therefore, at the critical level, one has
\be
\frac{d}{dc} \cR(x)\bigg\vert_{cr} = \cM(x) \frac{d}{dc} 
\frac{1}{\mu(xq^c)\mu(xq^{-c-2})\mu(x^{-1})^2} + \frac{1}
{\mu(xq^c)\mu(xq^{-c-2})\mu(x^{-1})^2} \frac{d}{dc} \cM(x) = 0 \,. 
\label{eq419}
\ee
Finally, the derivative of the term $T(x)$ is given by (using the same 
trick (\ref{eq415}))
\be
\frac{dT}{dc}(x)\bigg\vert_{cr} = -(\ln q) \left(x^{-1}\frac{d}{dx^{-1}}
\ln\tau(x^{-1}q^{1/2}) - x\frac{d}{dx}\ln\tau(xq^{1/2})\right) \,.
\ee
Now, from eqs. (\ref{eq412}) and (\ref{eq419}), 
$\displaystyle{\frac{d\cY}{dc}(x)}$ at the critical 
level is proportional to the identity matrix. The formula (\ref{eq41}) of
Theorem \ref{thmtwo} then immediately follows taking $x=z/w$.
\hfill \rule{5pt}{5pt}

\medskip

The structure of the Poisson bracket therefore derives wholly from the 
prefactor in the $\cY$ matrix (\ref{eq46}), a fact to be kept in mind. 
As a consequence, any dependence in $p$ is absent in the Poisson bracket 
structure function.

Note also that although we have not (and will not here) shown that the 
generators $t(z)$ generate the whole center at $c=-2$, we have nonetheless 
proved that their exchange algebra did generate a closed Poisson algebra 
for the corresponding set of classical variables.

\sect{Explicit Poisson structures}

{From} the equation (\ref{eq41}) and the periodicity property of the function 
$\tau(x)$, one gets easily
\bea
\Big\{ t(z),t(w) \Big\} &=& -(2\ln q) \left[ \sum_{n \ge 0} ~ \left( 
\frac{2x^2q^{4n+2}}{1-x^2q^{4n+2}} - \frac{2x^{-2}q^{4n+2}}{1-x^{-2}q^{4n+2}} 
\right ) \right. \nonumber \\
&& \left. +\sum_{n>0} ~ \left( 
- \frac{2x^2q^{4n}}{1-x^2q^{4n}} + \frac{2x^{-2}q^{4n}}{1-x^{-2}q^{4n}}\right)
- \frac{x^2}{1-x^2} + \frac{x^{-2}}{1-x^{-2}} \right] ~ 
t(z) \, t(w) \, , \label{eq51}
\eea
where $x=z/w$.
Interpretation of the formula (\ref{eq51}) must now be given in terms of 
the modes of $t(z)$, defined in the sense of generating function:
\be
t_n = \oint_C \frac{dz}{2\pi iz} \, z^{-n} \, t(z) \, .
\ee
The structure function $f(z/w)$ which defines the Poisson bracket 
(\ref{eq51}) has singularities which require a careful analysis of the 
procedure to get Poisson structures for $\{t_n, n \in \ZZ\}$. It is periodic 
with period $q^2$ and has simple poles at $z/w = \pm q^k$ for $k \in \ZZ$. 
In particular it is singular at $z/w = \pm 1$.

As a consequence, the expected definition of the Poisson structure 
$\{t_n,t_m\}$ as a double contour integral of (\ref{eq51}) must be made
more precise. Deformation of, say, the $w$-contour while the $z$-contour is 
kept fixed may induce the crossing of singularities of $f(z/w)$ which in turn 
modifies the computed value of the Poisson bracket. One can equivalently say 
that the explicit evaluation of $\displaystyle{\oint \frac{dz}{2\pi iz} \oint
\frac{dw}{2\pi iw} \, f(z/w) \, t(z) \, t(w)}$
requires an expansion of $f(z/w)$ in series of $z/w$, the form of which
in turn depends explicitly upon the value of $|z/w|$ being in a particular
interval $[q^k,q^{k+1}]$ for some $k \in \ZZ$.

Moreover, the singularity at $z/w= \pm 1$ implies that one cannot identify
$\displaystyle{\oint_{C_1} \frac{dz}{2\pi iz} \oint_{C_2}
\frac{dw}{2\pi iw}}$ with the permuted double integral
$\displaystyle{\oint_{C_2} \frac{dz}{2\pi iz} \oint_{C_1}
\frac{dw}{2\pi iw}}$ since by deformation of $C_1 \rightarrow C_2$ and 
$C_2 \rightarrow C_1$, the contours necessarily cross at some point. 
The immediate consequence is that the quantity
$\displaystyle{\oint_{C_1} \frac{dz}{2\pi iz} \oint_{C_2}
\frac{dw}{2\pi iw} \, z^{-n} \, w^{-m} \, f(z/w) \, t(z) \, t(w)}$
is {\em not} antisymmetric under the exchange $n \leftrightarrow m$
and cannot be taken as a Poisson bracket.

This leads us to define the Poisson bracket as:
\begin{defin}\label{defone}
\be
\{ t_n,t_m \} = \frac{1}{2} \left( \oint_{C_1} \frac{dz}{2\pi iz} \oint_{C_2} 
\frac{dw}{2\pi iw} + \oint_{C_2} \frac{dz}{2\pi iz} \oint_{C_1} 
\frac{dw}{2\pi iw} \right) z^{-n} \, w^{-m} \, f(z/w) \, t(z) \, t(w) \, .
\label{eq53}
\ee
\end{defin}
Such a procedure guarantees the antisymmetry of the postulated Poisson
structure due to the property $f(z/w) = -f(w/z)$. 
Remark that this definition is a semi-classical limit of a 
well-known procedure in two-dimensional field theory by which one computes
commutators of field modes out of radial-ordered contour integrals of
their operator product expansions. Here the relevant exchange algebra
is given by (\ref{eq44})-(\ref{eq45}).

The problem of contour deformation crossing the singularities at 
$z/w = \pm q^k$ where $k \ne 0$ and modifying the Poisson bracket is solved 
as follows. To simplify the notations we now choose the contours $C_1$ and 
$C_2$ to be circles of radii $R_1$ and $R_2$ respectively. We now state:
\begin{prop}\label{proptwo}
For any $k \in \ZZ^+$ such that $R_1/R_2 \in [q^{\pm k},q^{\pm(k+1)}]$, 
Definition \ref{defone} defines a consistent Poisson bracket whose specific 
form depends on $k$.
\end{prop}
By ``consistent'' we mean ``antisymmetric and obeying the Jacobi identity''.
\\
To prove Proposition \ref{proptwo} we now evaluate the postulated Poisson
brackets.

\medskip

$\bullet$ case $k = 0$: The first four terms in (\ref{eq51}) can be expanded 
in convergent power series of $x^{\pm s}$ without further ado since they have 
no singularity in this region and $x \in [q\, , \, q^{-1}]$ ($|q|<1$), hence 
the expansion parameters $x^{\pm 2}q^{4n+2}$ ($n\geq 0$), $x^{\pm 2}q^{4n}$ 
($n>0$) are smaller than $1$. Their contribution to the Poisson structure is
\be
\{ t_n,t_m \}_{k=0} = -4\ln q \oint_{C_1} \frac{dz}{2\pi iz} \oint_{C_2} 
\frac{dw}{2\pi iw} \sum_{s > 0} \left ( \frac{q^{4s}-q^{2s}}{1-q^{4s}}
x^{-2s} + \frac{q^{2s}-q^{4s}}{1-q^{4s}} x^{2s} \right) 
z^{-n} \, w^{-m} \, t(z) \, t(w) \,, \label{eq54}
\ee
where $C_1,C_2$ can be arbitrarily deformed once the contour 
integrals are evaluated by residue method. The last two terms have poles 
at $z/w = \pm 1$ and require two distinct expansions depending whether 
$R_1 > R_2$ or $R_2 > R_1$. Their contribution can then be 
reintroduced in the formal series expansion (\ref{eq54}) to give:
\be
\{ t_n,t_m \}_{k=0} = -2 \ln q \oint_{C_1} \frac{dz}{2\pi iz} \oint_{C_2} 
\frac{dw}{2\pi iw} \sum_{s \in \ZZ} \frac{q^s-q^{-s}}{q^s+q^{-s}}
\Big( \frac{z}{w} \Big)^{2s} z^{-n} \, w^{-m} \, t(w) \, t(z) \, ,\label{eq55}
\ee
which, evaluated by residues, gives an unambiguous answer. Remarkably
this is the (not centrally extended) Poisson bracket structure found
in the trigonometric (albeit with elliptic prefactor) case using a
completely different approach \cite{FR}. This identification requires 
the use of the parity of $t(z)$ to redefine modes as
$\displaystyle{\tilde t_n \equiv \oint \frac{dz}{2\pi iz} \, z^{2n} \, t(z)}$ 
in which case (\ref{eq55}) is exactly formula (9.4) in \cite{FR}. It is indeed 
a consistent Poisson structure and its identification establishes a 
connection between the elliptic and trigonometric $q$-deformed algebras, to 
be commented at the end.

\medskip

$\bullet$ case $k \ne 0$:
The Poisson bracket (\ref{eq53}) can be evaluated as the sum of the $k=0$ 
bracket (\ref{eq55}) and the contributions to (\ref{eq53}) of the poles
at $z/w = \pm q^{\pm n}$ where $n = 1,\dots,k$, from the deformation of the
contour $C_1$ (for instance, keeping $C_2$ fixed). This contribution can be 
evaluated step by step. The deformation of $C_1$ from sector $k-1$ to sector 
$k$ adds to (\ref{eq55}) the contribution from the poles 
$z/w =\pm q^{\pm k}$.
This contribution is given by a contour integral around each pole, with
an overall sign depending on the parity of $k$. It reads as:
\be
-(-1)^k \ln q\oint_C \frac{dz}{2\pi iz} \, 2t(z) \, t(zq^k) \, z^{-n-m} \, 
q^{-km} +(-1)^k \ln q \oint_C \frac{dz}{2\pi iz} \, 2t(z) \, t(zq^{-k}) \, 
z^{-n-m} \, q^{km}
\ee
or equivalently, defining the distribution $\delta(u) = \sum_{n \in \ZZ}u^n$:
\be
-(-1)^k \ln q \oint_{C_1} \frac{dz}{2\pi iz} \oint_{C_2} \frac{dw}{2\pi iw}
\, t(z) \, t(w) \, z^{-n} \, w^{-m} \left( \delta\Big(\frac{w}{zq^k}\Big) 
- \delta\Big(\frac{wq^k}{z}\Big) + \delta\Big(-\frac{w}{zq^k}\Big) 
- \delta\Big(-\frac{wq^k}{z}\Big) \right) \, , \label{eq57}
\ee
where, as in (\ref{eq55}), the integrals are evaluated term by term using
residue formula and therefore the relative position of the contours $C_1$
and $C_2$ is not relevant. Again, factorization of $\delta(zq^k/w) +
\delta(-zq^k/w)$
is made possible by the parity properties of $t(z) \equiv t(-z)$.
Interestingly enough, when one plugs back (\ref{eq57}) into (\ref{eq55})
using the formal series definition for the $\delta$ distribution, one gets
a compact formula very similar to (\ref{eq55}):
\be
\{ t_n,t_m \}_k = (-1)^{k+1} 2 \ln q \oint_{C_1} ~ \frac{dz}{2\pi iz} \oint_{C_2}~ \frac{dw}{2\pi iw} \sum_{s \in \ZZ} \frac{q^{(2k+1)s} - q^{-(2k+1)s}}
{q^s+q^{-s}} \Big( \frac{z}{w} \Big)^{2s} z^{-n} \, w^{-m} \, t(w) \, t(z) \, ,
\label{eq58}
\ee
which we now take as the final generic definition of the Poisson brackets
deduced from (\ref{eq51}). One indeed shows easily that for all $k$,
(\ref{eq58}) obeys Jacobi identity and is antisymmetric.

\medskip

\sect{Conclusion}

We have obtained a family, indexed by an integer $k$, of consistent Poisson 
brackets defined on the center $\cZ$ at $c=-2$ of $\elp$. The form of these 
Poisson brackets is similar to the form of Poisson bracket obtained by 
\cite{FR}.

We would like to expand a little more on this similarity. It is
difficult to take a direct limit $p \rightarrow 0$ of our algebra to get the
trigonometric $q$-deformed algebra, since this limit is highly singular 
\cite{JMK}. This seems to preclude a direct derivation from our algebra 
to the one obtained in \cite{FR}.

However we think that fruitful comparison can already be made at this point.

The occurrence of an elliptic Poisson algebra in \cite{FR} starting from a 
trigonometric $R$-matrix is probably connected to the introduction of an 
explicit elliptic prefactor in the initially trigonometric $R$-matrix. 
Our computation has indeed shown that the elliptic Poisson brackets entirely 
arise from the $\tau$-prefactor in our $R$-matrix. A similar phenomenon might 
then occur in \cite{FR}.

Two apparent discrepancies between the two structures can also be understood. 
Frenkel and Reshetikhin obtain one single Poisson bracket structure, which 
seems to correspond more closely to our Poisson bracket at $k=1$, but with 
a purely central $\delta$-term instead of the $t(z)t(zq)$ term in 
(\ref{eq57}). Remember however that the derivation in \cite{FR} uses an 
explicit representation of the trigonometric $q$-deformed algebra in terms 
of free quasi-bosons. It is possible that this particular representation, 
combined with the non trivial procedure leading from full elliptic to 
trigonometric algebra, leads to a degeneracy of such terms as $t(z)t(zq)$ 
and contains implicitely the restriction to the $k=1$ bracket. Provided 
that the limit procedure $p\rightarrow 0$ be better understood, this 
indicates a scheme through which the results in \cite{FR} could be 
connected to our general setting.  

\medskip

{\bf Acknowledgements.}
This work was supported in part by CNRS and Foundation Angelo della Riccia. J.A.
wishes to thank ENSLAPP-Lyon and CRM-Montr\'eal for their hospitality. 
L.F. is indebted to Centre de Recherches Math\'ematiques of Universit\'e
de Montr\'eal for its kind invitation and support. M.R. thanks ENSLAPP-Lyon 
for its kind hospitality. 
We have benefited from fruitful discussions with L. Freidel and J.M. Maillet.

\end{document}